\documentclass[11pt,a4paper]{article}
\usepackage{fouriernc}
\usepackage{tabularx}
\usepackage{booktabs}
\usepackage{makecell}
\usepackage{amsmath, amssymb}
\usepackage{graphicx}
\usepackage{xcolor,ulem}
\usepackage{siunitx}
\usepackage[export]{adjustbox}
\usepackage{subfig}
\usepackage[utf8]{inputenc}
\usepackage[T1]{fontenc}
\usepackage{enumitem}
\usepackage{color, soul}
\usepackage{authblk}
\usepackage{marvosym}
\usepackage{wasysym}
\usepackage{multirow}
\usepackage{url}
\usepackage{apalike}
\usepackage{hyperref}
\hypersetup{colorlinks=true,allcolors=blue}
\usepackage{vmargin}
\setmarginsrb{2cm}{2cm}{2cm}{2cm}{0mm}{0mm}{0mm}{10mm}
\usepackage{tikz}
\usepackage{pgfplots, pgfplotstable}
\pgfplotsset{compat=1.7}

\title{\vspace*{-1.5cm} \bfseries Enhancing academic performance:\\ The impact of active learning in mathematical economics}
\author[1]{P. K. Ng\,}
\author[2]{N. Karjanto\thanks{\Letter: \url{natanael@skku.edu}\href{https://orcid.org/0000-0002-6859-447X}{\includegraphics[scale=0.08]{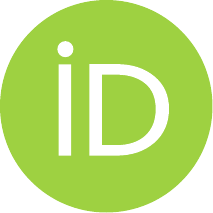}}}}
\affil[1]{Centre for English Language and Foundation Education, The University of Nottingham Malaysia \protect \\ Jalan Broga, Semenyih 43500, Selangor, Malaysia}
\affil[2]{Department of Mathematics, University College, Natural Science Campus, Sungkyunkwan University\protect \\ 2066~Seobu-ro, Jangan-gu, Suwon~16419, Gyeonggi-do, Republic of Korea}

\date{\vspace*{-0.5cm} \footnotesize Last updated \today}

\begin{document}
\maketitle
\vspace*{-0.5cm}
\begin{abstract}
This paper explores the impact of active learning in mathematical economics on students' academic performance (assessment scores). An experimental design involving foundation students enrolled in the arts and business and management foundation programmes in a British university located in Malaysia was adopted. The control group underwent the more traditional lecture method with the students taking on a passive role of listening to information disseminated by the instructor. The treatment group, in contrast, was given minimum explanation with the bulk of learning coming from students actively solving problems presented in case studies based on real-world events. Results show that the 189 students in the treatment group performed significantly better than the 146 students in the control group.\\

\noindent
Keywords: active learning, academic performance, mathematical economics, university foundation, British university in Malaysia.
\end{abstract}

\section{Introduction}
As higher education institutions gear towards producing more well-balanced graduates who are able to think both logically as well as creatively and have a more holistic/broader worldview, students are increasingly given more flexibility in choosing modules from different fields. It is thus, not uncommon at all to see mathematical economics being offered to non-economics major students.
 
One prevalent grouse that economics students have, including those specialising in the subject, is that the constructs/concepts are often too theoretical or abstract to understand, and they do not see the value of these models in terms of their relevance to the real world. Consequently, they become uninterested in, or disengaged with, the module. The major challenge that teachers have, therefore, is to design and deliver engaging pedagogy that ideally, could foster sustainable learning practices that extend beyond the completion of the module or graduation from the university--a task that is made all the more arduous when the module is offered to students with different levels of skills, knowledge, ability and interests in mathematics and economics. In this regard, active learning is touted as a pedagogy with the potential to revolutionise the teaching and learning environment, paving the way for more inquisitive and creative thinkers. Indeed, extensive research literature consistently affirms the effectiveness of active learning pedagogy, particularly in the context of mathematical economics education, in enhancing academic performance, improving information retention, and developing critical thinking skills. See e.g.,~\cite{hettler2015active,ray2018teaching,sekwena2023active}.

Active learning pedagogy is an educational approach that places the learner at the center of the learning process, encouraging them to actively engage in, and take responsibility for, the learning experience rather than being passive recipients of information~\cite{anthony1996active,gogus2012active}. It involves a shift away from traditional lecture-style teaching, where the instructor delivers information to the students, towards a more interactive and participatory learning environment. In active learning, learners are encouraged to think critically, solve problems, and construct their knowledge actively through various activities and discussions. The instructor plays the role of a facilitator or guide, helping students explore concepts, ask questions, and develop a deeper understanding of the subject matter. Key features of active learning pedagogy include, but are not limited to, a student-centered approach, engagement and participation, collaborative learning, critical thinking and problem-solving, reflection and feedback, use of technology, continuous assessment, inclusive learning environment, and application of knowledge~\cite{freeman2014active,karjanto2019english,karjanto2022sustainable,kerrigan2018active}.

This research paper seeks to delve into the profound influence of active learning strategies in mathematical economics education, with the central aim of investigating the following research question, ``What is the comparative impact on academic performance between the experimental group that embraces active learning activities and the control group that receives traditional lectures?''

\section{Research methodology}

In order to investigate the effect of an active learning instructional approach on students' academic performance, we used the independent measures (between-groups) experimental design involving foundation students majoring in arts and business and management at a private university in Malaysia--an overseas campus branch of a reputable British university.  Both the control and treatment groups attended a 1.5-hour lecture session twice a week for 5~weeks in Summer 2022 and 2023 respectively. Students in the control group were taught the more traditional way, i.e., mainly instructor-led sessions with detailed explanations of theories/models (e.g., demand-supply, production possibility frontier, and aggregate demand-aggregate supply) as represented by graphs. Short case studies of events from the real world were then used by the instructor to illustrate how these models could be used in reality. The treatment group, in contrast, was given a very brief introduction to the economic models with the hope that most of the learning and understanding would come from students solving problems by applying economic models to the case studies drawn from real-world events. The students are given time to work through these exercises on their own first, or they could discuss them with the person seated next to them, with the instructor moving around the classroom to check on their progress and clarify doubts, before finally discussing the questions with the whole class. 

At the end of the course, both the control and treatment groups were assessed through a group assignment and an online test, with both components added to give the students an overall score.

\section{Results}

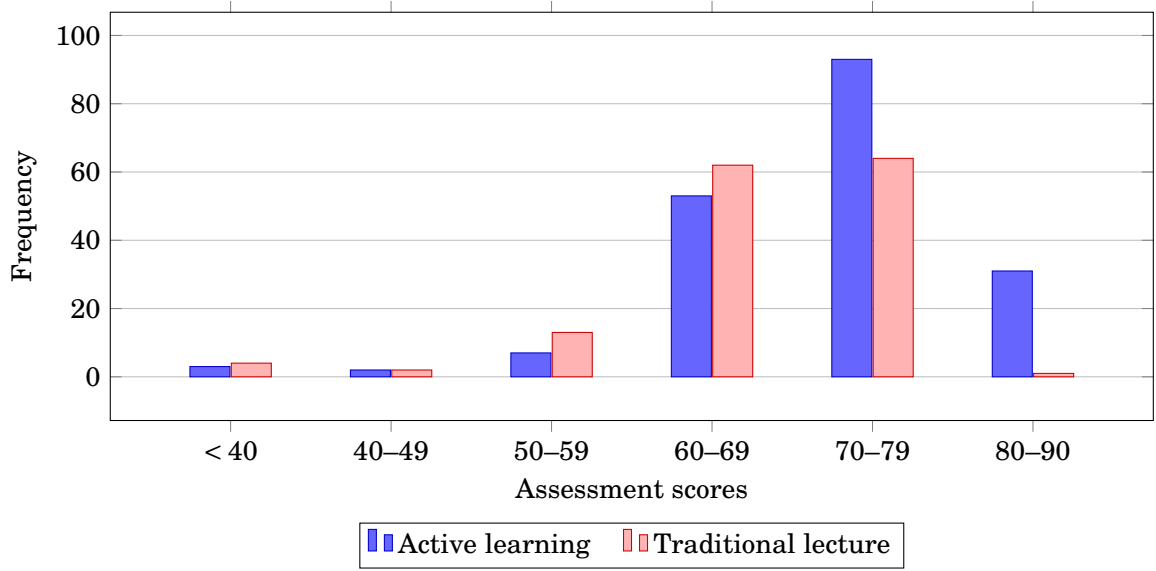
\begin{figure}[htbp]
\centering
\begin{tikzpicture}
\begin{axis}[width=0.9\textwidth,height=7cm,	
			ybar=0.5pt,
			bar width=15pt,
			enlargelimits=0.15,
			xlabel={Assessment scores},
			symbolic x coords={$<40$,40--49,50--59,60--69,70--79,80--90},
			xtick=data,
			ylabel={Frequency},
			ymajorgrids,
			legend style={
				at={(0.5,-0.25)},
				anchor=north,legend columns=-1}
			]
\addplot[fill=blue!60!white,draw=blue!80!black] coordinates {($<40$,3) (40--49,2) (50--59,7)  (60--69,53) (70--79,93) (80--90,31)};
\addplot[fill= red!30!white,draw=red!80!black ] coordinates {($<40$,4) (40--49,2) (50--59,13) (60--69,62) (70--79,64) (80--90,1)};
\legend{Active learning{$\quad$}, Traditional lecture}
\end{axis}
\end{tikzpicture}
\caption{Frequency of students according to assessment scores in the active learning group vs. the traditional lecture group.} 	\label{figurehistogram}
\label{figure}
\end{figure}

We used the two-sample t-test to determine whether there is any significant difference in assessment scores between those who have undergone the active learning method vs. those who were lectured the more traditional way. The frequency of students according to assessment scores is illustrated in Figure~\ref{figurehistogram}, while the outputs of the statistical tests are presented in Table~\ref{tablegroupstatistics} (showing group statistics) and Table~\ref{tablettest} (results of the independent samples t-test).
\begin{table}[h]
\centering
\begin{tabular}{@{}cccccc@{}}
\toprule	
						& Group					 & $N$ & Mean & Standard deviation & Standard error \\ \midrule
Academic performance 	& Active learning		 & 189 & 71.14427 & 8.285017 & 0.602646 \\	\cline{2-6}
(Assessment scores) 	& Traditional lecture	 & 146 & 66.63959 & 9.362886 & 0.774878 \\				
\bottomrule
\end{tabular}
\caption{Group statistics.}			\label{tablegroupstatistics}
\end{table}
\begin{table}[h]
\centering
\begin{tabular}{@{}ccccccc@{}}
\toprule	
	$t$	& df  & $p$      & \makecell[c]{Mean\\ difference} & \makecell[c]{Standard\\ error\\ difference} & \multicolumn{2}{c}{\makecell[c]{95\% confidence interval\\ of the difference}} \\ \cline{6-7}
		&	  &		     &				   &						   & Lower 	& Upper \\ \midrule
4.6616	& 333 & 0.000005 & 4.505           & 0.966                     & 2.6037 & 6.4057 \\				
\bottomrule
\end{tabular}
\caption{Results from two sample $t$-test.}			\label{tablettest}
\end{table}

Results show that the 189 students who experienced active learning (mean $= 71.14$, SD $= 8.28$) compared to 146 students in the control group (mean $= 66.64$, SD $= 9.36$) demonstrated a significantly better academic performance $t\left(333\right) = 4.66$, $p < 0.00001$.

\section{Discussion and conclusion}

The findings in this study demonstrate that active learning via case studies that allows students to discuss and actively apply economic models (i.e., graphical analysis) to understand real-world events is more effective in boosting academic performance compared to the more traditional teacher-dominated pedagogy. Thus, it adds to the growing literature in support of anchoring teaching practices on active learning strategies. See e.g.,~\cite{grabinger1995rich,kerrigan2018active,ray2018teaching,yadav2014case}.  

Ideally, in this type of case-based method, the class should be split into smaller groups of 3--5 students with all groups being able to present/share their answers~\cite{brame2016active}. However, owing to the relatively large class size of about 100 students and time constraints, this would not have been practical. In our modified version, students are given the flexibility to either work on the cases on their own or discuss them with the person next to them. Several challenges emerged from this, notably:
\begin{enumerate}[leftmargin=1.4em]
\item Some students were not participating, i.e., they did not appear interested in solving the questions.
\item The instructor could not attend to/check on all the students despite some of them requiring further attention/guidance.
\item Some students found the whole exercise confusing as they were unsure of where/how to start and they did not quite trust the information/views of their peers. They expressed a preference for having the instructor lead the session; guiding them step-by-step on how to apply the economic models.
\end{enumerate}

However, these issues should not be seen as a deterrent to implementing active learning in relatively large classes. Instead, more careful planning is needed for future implementations. For instance, informing students in advance to form groups so that they can better organise their seating arrangements during class. In addition, less important topics could be removed so that the instructor has more time to check on the understanding of students. A questionnaire could be given to students at the end of the course to better gauge their subjective assessment of this pedagogy. Also, given the diverse backgrounds of the students, some of whom have never taken economics in high school, a diversified teaching approach is needed, perhaps with the instructor taking on a more active role in guiding those who are new to the subject area. 

\subsection*{Conflict of Interest}
The authors declare that they have no conflicts of interest.

\end{document}